\documentclass[11pt]{amsart}
\usepackage{amsmath}
\usepackage{latexsym}
\usepackage{amssymb}
\usepackage{amscd}
\input epsf

\usepackage{amsfonts}
\usepackage[all]{xy}

\newtheorem{proposition}{Proposition}[section]

\newtheorem{lemma}[proposition]{Lemma}
\newtheorem{definition}[proposition]{Definition}
\newtheorem{lemma-definition}[proposition]{Lemma-Definition}
\newtheorem{theorem}[proposition]{Theorem}

\newtheorem{corollary}[proposition]{Corollary}

\newtheorem{example}[proposition]{Example}
\newtheorem{remark}[proposition]{Remark}

\newenvironment{dok}{\par\vspace{-5pt}%
\par\noindent\begingroup%
\leftskip=0em\hspace{0em}{\bf Proof.}}%
{\endgroup\hfill$\Box$}

\newcounter{tmp}

\def\db#1{ \bD^b({#1})}

\def\perf#1{{\mathfrak P}{\mathfrak e}{\mathfrak r}{\mathfrak f}({#1})}
\def\dsing#1{ \bD_{\rm Sg}({#1})}
\def\dsingq#1{ \bD^{'}_{\rm Sg}({#1})}

\def\ds#1{ {#1}_{\rm Sg}}
\def\dhf#1{ {#1}_{\rm hf}}

\def\ove{\overline}

\def\Hom{{\mathrm H}{\mathrm o}{\mathrm m}}

\def\h#1,#2{{\operatorname{Hom}}({#1}\:,\; {#2})}
\def\Ho#1,#2,#3,#4{{\operatorname{Hom}}^{#1}_{#2}({#3}\:,\; {#4})}
\def\Ex#1,#2,#3,#4{{\operatorname{Ext}}^{#1}_{#2}({#3}\:,\; {#4})}
\def\Mod{\operatorname{Mod}\!}
\def\mod{\operatorname{mod}\!}
\def\Proj{\operatorname{Proj}\!}
\def\proj{\operatorname{proj}\!}
\def\coh{\operatorname{coh}}
\def\Qcoh{\operatorname{Qcoh}}

\def\lto{\longrightarrow}

\def\C{{\mathcal C}}

\def\D{{\mathcal D}}
\def\F{{\mathcal F}}
\def\E{{\mathcal E}}
\def\G{{\mathcal G}}

\def\N{{\mathcal N}}
\def\O{{\mathcal O}}

\def\L{{\mathcal L}}

\def\P{{\mathcal{P}}}

\def\ZZ{{\mathbb Z}}

\def\bD{{\mathbf D}}
\def\bR{{\mathbf R}}

\def\bL{{\mathbf L}}

\def\ZZ{{\mathbb Z}}
\def\b1{{\mathbf 1}}

\def\AA{{\mathbb A}}

\def\ZZ{{\mathbb Z}}

\def\PP{{\mathbb P}}

\def\Hom{\operatorname{Hom}}

\def\Ext{\operatorname{Ext}}

\def\ext{\underline{{\mathcal E}xt}}
\def\tor{\underline{{\mathcal T}or}}
\def\rhom{\bR\underline{{\mathcal H}om}}

\def\Sing{\operatorname{Sing}}

\def\Spec{\operatorname{Spec}}

\def\id{{\operatorname{id}}}

\title[]{Triangulated categories of singularities and
equivalences between Landau-Ginzburg models}

\author[]{Dmitri Orlov}

\address{ Algebra Section, Steklov Mathematical Institute RAS,
Gubkin str. 8, Moscow 119991, RUSSIA}

\email{orlov@mi.ras.ru}

\thanks{This work is done
under partial financial support of the Weyl Fund,
the grant RFFI (No~02-01-00468), grant of President of RF
in support of young russian scientists MD-2731.2004.1,
grant CRDF Award No RM1-2405-MO-02. It is also a
pleasure for me to express my gratitude to the Russian Science
Support Foundation.}

\date{}

\begin{document}
\begin{abstract}
In this paper we prove an existence of some type of equivalences  between
triangulated categories of singularities for varieties
of different dimensions. This class of equivalences
generalizes so called Kn\"orrer periodicity.
As consequence we get  equivalences between categories
of D-branes of type B on Landau-Ginzburg models of
different dimensions.
\end{abstract}
\maketitle
\section*{Introduction}

The bounded derived
category of coherent sheaves $\db{\coh(X)}$
is a natural triangulated category
which can be attached to  an algebraic variety $X.$
This category
has a triangulated subcategory $\perf{X}$ formed by  perfect
complexes.
A notion of a perfect complex  was introduced in \cite{Il},
and, by definition, it is a complex of  sheaves of $\O_X$\!-modules which locally is
quasi-isomorphic to a bounded complex of locally free sheaves of finite
type (a good reference is \cite{TT}).

If the variety $X$ is smooth then any coherent sheaf
has a finite resolution of locally free sheaves of finite type and
the subcategory of perfect complexes coincides with
$\db{\coh(X)}.$ But for singular varieties this property is not
fulfilled. One can define  a triangulated category
of singularities $\dsing{X}$ as the quotient of the triangulated
category $\db{\coh(X)}$ by the full triangulated subcategory of
perfect complexes $\perf{X}$ \cite{Tr}. The category $\dsing{X}$ reflects
the properties of the singularities of $X$ and "does not depend on
all of $X$". For example it is invariant with
respect to a localization in  Zariski topology (\cite{Tr}, and Prop.
\ref{locpr}).

The investigation  of  such  categories  is not only connected with
a study of singularities but is mainly  inspired by the
Homological Mirror Symmetry Conjecture \cite{K}.

The HMSC has dealings with  Calabi-Yau  varieties
endowed with symplectic forms.
It asserts that if two Calabi-Yau varieties $X$ and
$Y$ are mirror symmetric to each other, then the  category of
D-branes of type B on $X$ is equivalent to the category of
D-branes of  type A on
$Y$, and vice versa.
From the mathematical point of view the category of D-branes of type
B on a variety is the derived category of coherent sheaves on it
\cite{K, Do}. As a candidate for a category of A-branes on
a symplectic manifolds the so-called Fukaya category has been proposed.
Its objects are, roughly speaking, Lagrangian submanifolds
equipped with flat vector bundles \cite{K}.

On the other hand, physicists also consider
 mirror symmetry relation for other varieties
 for instance Fano varieties.
 In these cases  the so called  Landau-Ginzburg theories arise
on the  mirror side
\cite{HV}. General definition of a Landau-Ginzburg model
involves, besides a choice of a target space, a choice of a
holomorphic function $W$ on it which is called superpotential.

For Fano varieties one has the derived categories of coherent
sheaves (B-branes) and given  a symplectic form, one
can propose a suitable Fukaya category (A-branes).
Thus, if one wants to extend the Homological Mirror  Symmetry
Conjecture to the non-Calabi-Yau case, one should  understand
D-branes in Landau-Ginzburg models.

Categories of A-branes in Landau-Ginzburg models are studied in
\cite{HIV} and in \cite{Se} from the mathematical point of view.
Mirror symmetry relates B-branes on a Fano variety (coherent
sheaves) to A-branes in a LG model.
One can also consider the  Fukaya category (A-branes) on a Fano
variety and can expect  in this case that the Fukaya category is
equivalent to the category of B-branes in the mirror
Landau-Ginzburg model.

A mathematical definition for the category of B-branes in
Landau-Ginzburg models was proposed by M.Kontsevich. Roughly, he
suggests that the superpotential $W$ deforms complexes of coherent
sheaves to "twisted" complexes, i.e the composition of
differentials  is no longer zero, but  is equal to multiplication
by $W$.
In the paper
\cite{Tr} we established a connection between  categories of B-branes in
Landau-Ginzburg models and triangulated categories of
singularities. We considered  singular fibres of the map $W$ and
showed that the triangulated categories of singularities of these
fibres are equivalent to the categories of B-branes.
In particular this gives
that the category of B-branes  depends only on the singular fibres of the
superpotential.

In this paper we  establish equivalence between
Landau-Ginzburg models
of different dimensions. More precisely,
we prove an equivalence between triangulated categories of singularities
of two different schemes. The first one is the zero subscheme of
a regular section $s\in H^0(S, \E)$ of a vector bundle $\E$
on a smooth scheme $S$ and the second one is
the respective divisor on the projective bundle $\PP(\E^{\vee})$
(Theorem
\ref{main1}).
This  result implies an equivalence of categories
of D-branes of type B for different Landau-Ginzburg models.
In particular,  it says the following: Let $W: Y\to\AA^1$ be a
superpotential on a variety $Y=S\times \AA^1$
of the form $W= f+xg,$ where $f,g$ are
 functions on $S$ and $x$ is an
coordinate on $\AA^1$.
Then one can reduce dimension and change
the Landau-Ginzburg model on $Y$ by the Landau-Ginzburg
model on $X\subset S,$ where $X$ is the zero subvariety
of the function $g,$ with a restriction of $f$ as a superpotential.
Thus, we prove that the category of B-branes on the
Landau-Ginzburg model on $Y$ is equivalent to the category of
B-branes in the Landau-Ginzburg model on $X$
(Corollary \ref{LG}).
Of course, we expect that the categories of
D-branes of type A are also equivalent for  these LG-modles.

I am grateful to  Anton Kapustin and Ludmil Katzarkov for
useful discussions. I also thank to Valery Lunts for  reading of a
preliminary draft of the paper and making a number of valuable
comments.
I also would like to express my gratitude to
Alexander Kuznetsov who pointed out an existence
of semiorthognal decomposition which arises in Proposition
\ref{dual}. I want to thank
Institute for Advanced Study, where this work was written,
for  hospitality and  very stimulating atmosphere.

\section{Triangulated categories of singularities}

We remind that a triangulated category
$\D$  is  an additive category with  the following data:
\begin{list}{\alph{tmp})}%
{\usecounter{tmp}}
\item an additive autoequivalence
$[1]: \D\lto\D$ (it is called a translation functor),
\item a class of exact (or distinguished) triangles:
$$
X\stackrel{u}{\lto}Y\stackrel{v}{\lto}Z\stackrel{w}{\lto}X[1],
$$
\end{list}
which must satisfy a certain set of axioms (see \cite{Ve}, also \cite{GM, KS, Ke, Nee}).

A functor $F : {\D} \lto{\D}'$ between two triangulated categories
${\D}$ and ${\D}'$ is called {\sf exact} if it commutes with the
translation functors and
transforms  exact triangles into exact
triangles.

Let $\N\subset \D$
be a full triangulated subcategory, i.e. it is a full subcategory
which closed with respect to the translation functor and it is triangulated
with respect to exact
triangles in $\D.$
Denote by $\Sigma(\N)$ a class of morphisms $s$ in $\D$
embedding into an exact triangle
$$
X\stackrel{s}{\lto} Y\lto N\lto X[1]
$$
with $N\in \N$.
It can be checked that $\Sigma(N)$ is a multiplicative system.
Define the quotient
$
\D/\N
$
as
a localization $\D[\Sigma(\N)^{-1}]$ \cite{GZ, GM, Ve}.
We endow the category $\D/\N$ with a translation functor
induced by the translation functor in the category $\D$.
The category $\D/\N$ becomes a triangulated category by taking for
exact triangles the images of exact
triangles in $\D$. The quotient functor $Q:\D\to \D/\N$
annihilates $\N$. Moreover, any exact functor $F: \D\to \D'$ between
triangulated categories for which $F(X)\simeq 0$ when $X\in \N$
factors uniquely through $Q$.

The following lemma, which will be necessary in the future, is evident.
\begin{lemma}\label{adjqu}
Let $\N$ and $\N'$ be full triangulated subcategories
of  triangulated categories $\D$ and $\D'$ respectively.
Let $F: \D\to \D'$ and $G: \D'\to \D$ be  adjoint pair
of exact functors such that $F(\N)\subset \N'$ and $G(\N')\subset \N$.
Then they induce functors
$$
\ove{F}:\D/\N\lto\D'/\N',
\qquad
\ove{G}:\D'/\N'\lto \D/\N
$$
which are adjoint too. Moreover, if the functor $F: \D\to \D'$ is
fully faithful then the functor $\ove{F}:\D/\N\lto\D'/\N'$ is fully faithful too.
\end{lemma}

We are mainly interested in triangulated categories and
their quotient by triangulated subcategories which are coming from algebraic
geometry.

Let $X$ be a scheme over field $k$. We say that $X$ satisfies the
condition (ELF) if it is

\vspace{0pt}
\begin{tabular}{ll}
(ELF)& \begin{tabular}{l}
 separated  noetherian of finite  Krull dimension and  has enough  locally free sheaves,\\
i.e. for any coherent sheaf $\F$
there is an epimorphism $\E\twoheadrightarrow\F$
for a locally free\\
 sheaf $\E.$
\end{tabular}
\end{tabular}

\vspace{3pt}
For example, any quasi-projective scheme satisfies these conditions.

Denote by $\db{\coh(X)}$ (resp. $\db{\Qcoh(X)}$)
 the bounded derived categories of coherent (resp. quasi-coherent)
 sheaves on $X$.
Since $X$ is noetherian the natural functor $\db{\coh(X)}\lto
\db{\Qcoh(X)}$ is fully faithful and realizes an equivalence of
$\db{\coh(X)}$ with the full subcategory
$\db{\Qcoh(X)}_{\coh}\subset \db{\Qcoh(X)}$ consisting of all
complexes with coherent cohomologies (see \cite{Il} II, 2.2.2).

The objects of the category
$\db{\coh(X)}$ which are isomorphic to bounded complexes of locally free sheaves
on $X$ form a full triangulated subcategory. It is called the subcategory
of perfect complexes and
denoted by $\perf{X}.$\footnote{
Actually, a perfect complex is defined as a complex of
$\O_X$\!-modules locally quasi-isomorphic to a bounded complex of
locally free sheaves of finite type. But under our assumption on
the scheme any such complex is quasi-isomorphic to a bounded
complex of locally free sheaves of finite type(see \cite{Il} II,
or \cite{TT} \S 2).}
\begin{definition}\label{trcsin} We define a triangulated category
$\dsing{X}$ as the quotient
of the triangulated category $\db{\coh(X)}$ by the full triangulated
subcategory
$\perf{X}$ and call it  a triangulated category of singularities of $X$.
\end{definition}
It is  known that if our scheme $X$  is
regular  then the
subcategory of perfect complexes coincides with the whole bounded
derived category of coherent sheaves. In this case the triangulated category
of singularities is trivial.

Let $f: X\to Y $ be a morphism
of finite Tor-dimension (for example a flat morphism or a regular closed embedding). It  defines
the inverse image functor $\bL f^*:\db{\coh(Y)}\lto \db{\coh(X)}$.
It is clear that the functor $\bL f^*$ sends
 perfect complexes on $Y$ to  perfect complexes on $X$.
Therefore, the functor $\bL f^*$ induces  an exact functor
$\bL \bar{f}^* :\dsing{Y} \lto \dsing{X}.$

Suppose, in addition, that the morphism  $f: X\to Y$ is  proper
and locally of finite type. Then  the functor
of direct image $\bR f_* : \db{\coh(X)}\to \db{\coh(Y)}$
takes  perfect complexes on $X$ to  perfect
complexes on $Y$ (see \cite{Il} III, or \cite{TT}).
Hence it determines a functor
$\bR \bar{f}_* :\dsing{X}\to \dsing{Y}$ which is the right
adjoint to $\bL \bar{f}^*.$

%We also  can consider  the full
%triangulated subcategory $\lfr{X}\subset \db{\Qcoh(X)}$ consisting of objects which are
%isomorphic to bounded complexes of locally free sheaves in
%$\db{\Qcoh(X)}$. Define a triangulated
%category $\dsingq{X}$ as
% the quotient $\db{\Qcoh(X)}/\lfr{X}$.
%The fully faithful functor $\db{\coh(X)}\lto \db{\Qcoh(X)}$
%induces a  functor $\dsing{X}\lto \dsingq{X}$
%which is also fully faithful (see \cite{Tr}).

A fundamental  property of  triangulated categories of singularities is a property of locality.
\begin{proposition}{\rm \cite{Tr}}\label{locpr}
Let X satisfy (ELF) and let $j: U\hookrightarrow X$ be an embedding
of an open subscheme such that $\Sing(X)\subset U$. Then the
functor $\bar{j}^*:\dsing{X}\lto \dsing{U}$ is an equivalence of
triangulated categories.
\end{proposition}

\begin{remark} {\rm The definition of triangulated categories of singularities
can be extended to orbifolds and more generally to stacks.
We briefly consider the example of quotient stacks.
Let $G$ be an affine group scheme of finite type which acts on a scheme
$S$ which satisfy (ELF). Consider the quotient stack $\left[ S/G \right].$
The category of coherent sheaves on this stack $\coh(\left[ S/G \right])$
 coincides
with the category $\coh^{G}(S)$ of $G$\!-equivariant coherent sheaves on
$S.$ The quotient of the bounded derived category of coherent sheaves
$\db{\coh(\left[S/G\right])}$
by the triangulated subcategory of perfect complexes
$\perf{\left[S/G\right]},$ which is formed by bounded complexes
of locally free sheaves in $\db{\coh(\left[S/G\right])}\cong
\db{\coh^G(S)},$ can be called
a triangulated category of singularities of the quotient stack
$\left[S/G\right].$ (For more information when the quotient stack
has enough locally free sheaves see papers \cite{Tom, Tot}.)}
\end{remark}
\begin{remark}{\rm
One also can consider noncommutative case.
Let $A$ be a right noetherian algebra. Denote by
$\mod-A$ and $\Mod-A$ the abelian categories
of finitely generated right modules and all right modules
respectively.
Consider bounded derived categories $\db{\mod-A}$ and
$\db{\Mod-A}.$ They have  triangulated subcategories
consisting of objects which are isomorphic to bounded complexes
of projectives. These subcategories can be considered as
 derived of exact categories of projective modules
 $\db{\proj-A}$ and $\db{\Proj-A}$ respectively (see for example \cite{Ke}).
Now we can define  triangulated categories of singularities $\dsing{A}$
and $\dsingq{A}$ as the quotients $\db{\mod-A}/\db{\proj-A}$
and $\db{\Mod-A}/\db{\Proj-A}$ respectively. As in the commutative case,
if $A$ has a finite homological dimension then
these quotient categories are trivial.}
\end{remark}

There is a  way to generalize  Definition \ref{trcsin}
to any triangulated category.
Let $\D$ be a triangulated category.
\begin{definition}We say that an object
$A$ is {\sf homologically finite} if for any object $B\in\D$
all $\Hom(A, B[i])$ are trivial except for finite number of $i\in \ZZ.$
All such objects form a full triangulated subcategory which
will be denoted by $\dhf{\D}.$
\end{definition}
\begin{definition}
We define a triangulated category
$\ds{\D}$ as the quotient $\D/\dhf{\D}$
of the triangulated category $\D$ by the full triangulated
subcategory
$\dhf{\D}.$
\end{definition}
These categories $\dhf{\D}$ and $\ds{\D}$ have  good behavior
with respect to semiorthogonal decomposition of $\D.$
We  recall some definitions and facts concerning
admissible subcategories and semiorthogonal decompositions (see \cite{BK, BO}).
\begin{definition}\label{adm}
Let $I:\N\hookrightarrow\D$ be an embedding of a full triangulated subcategory
$\N$ in a triangulated category $\D.$
We say that ${\N}$ is {\sf right admissible} (resp. {\sf left admissible})
if there is a right (resp. left) adjoint functor $P:\D\to \N.$
The subcategory $\N$  will be called admissible if it is right and left admissible.
\end{definition}

Let ${\N\subset\D}$ be a full triangulated subcategory.
The {\sf right orthogonal} to ${\N}$ is a full subcategory
${\N}^{\perp}\subset {\D}$ consisting of all objects $M$ such that ${\Hom(N, M)}=0$
for any $N\in{\N}$. The {\sf left orthogonal} ${}^{\perp}{\N}$ is defined analogously.
The orthogonals  are also triangulated subcategories.

The property to be right  admissible
for the subcategory $\N$ is equivalent to the following:
for each $X\in{\D}$ there is an exact triangle $N\to X\to M$,
where $N\in{\N}$ and $M\in{\N}^{\perp}.$

If $\N\subset\D$ is an admissible subcategory then we say that the category
$\D$  has  semiorthogonal decompositions of the form
$\left\langle{\N}^{\perp},{\N}\right\rangle$ and
$\left\langle{\N},{}^{\perp}{\N}\right\rangle.$
%Sometimes this process of decomposition can be extended on $\N$
%and its orthogonals.

\begin{definition}\label{sd}
A sequence of admissible subcategories $({\N}_1, \dots, {\N}_n)$ in a derived category ${\D}$
is said to be {\sf semiorthogonal} if the condition ${\N}_j\subset {\N}^{\perp}_i$ holds
when $j<i$ for any $1\le i\le n$.
In addition, a semiorthogonal sequence is said to be {\sf full} if it
generates the category ${\D}$.
In this case we call it as  a {\sf semiorthogonal decomposition}
of the category ${\D}$ and denote this as
$$
{\D}=\left\langle{\N}_1, \dots, {\N}_n\right\rangle.
$$
\end{definition}
\begin{proposition}
Suppose that a triangulated category ${\D}$ has
a semiorthogonal decomposition $\D=\left\langle{\N}_1,....,{\N}_n\right\rangle.$
Then the categories $\dhf{\D}$ and $\ds{\D}$  also have semiorthogonal decompositions
of the same form
\begin{equation}\label{decom}
\dhf{\D}=\left\langle\dhf{(\N_1)}, \dots, \dhf{(\N_n)}\right\rangle,
\qquad
\ds{\D}=\left\langle\ds{(\N_1)}, \dots, \ds{(\N_n)}\right\rangle.
\end{equation}
\end{proposition}
\begin{dok}
First, note that if a functor $u:\D\lto \D'$
has a right adjoint $v:\D'\lto D$ then
$u(\dhf{\D})\subset\dhf{\D'},$ because
$$
\Hom(uA, B[j])\cong \Hom(A, vB[j])
$$
for any $A\in\D,\; B\in \D'.$
Moreover, if $u$ is a full embedding then $u(\dhf{\D})=\dhf{\D'}\cap u(\D).$

Second, we can assume that $n=2.$
Since $\N_k$ are admissible the embedding functor $i_k: \N_k\lto\D$
has right adjoint $p_{k}$. Hence, it takes any homologically finite object
to a homologically finite object.
Now  suppose $X\in\dhf{\D}$ and consider decomposition
$$
i_1 p_1 (X)\lto X\lto i_2 q_2 (X),
$$
where $q_2:\D\lto \N_2$ is left adjoint to $i_2.$
As we proved above $q_2(X)\in\dhf{(\N_2)}.$
Hence, $i_2 q_2(X)\in\dhf{\D}.$ This implies that
$i_1 p_1(X)\in\dhf{\D}$ and, consequently, $p_1(X)\in \dhf{(\N_1)}.$
This gives us the semiorthogonal decomposition for $\dhf{\D}$ of the form
(\ref{decom}).
Finally, applying Lemma \ref{adjqu} we obtain the semiorthogonal decomposition for
$\ds{\D}$ of the form
$\ds{\D}=\left\langle\ds{(\N_1)}, \dots, \ds{(\N_n)}\right\rangle.$
\end{dok}

\begin{proposition}
Let $X$ satisfies (ELF). Then the subcategory
$\dhf{\D}$ of homologically finite objects in $\D\cong\db{\coh(X)}$
coincides with the subcategory $\perf{X}$ and, hence, $\ds{\D}\cong\dsing{X}.$
\end{proposition}
\begin{dok}
If an object $A\in\D$ is a perfect complex then it is quasi-isomorphic
to a bounded complex of vector bundles.
Since the cohomologies of any coherent sheaf is bounded by the
Krull dimension of variety we have that for any vector bundle $\P$ and any coherent sheaf $\F$
there is equality $\Ext^i(\P, \F)=0$ when $i$ is greater than Krull dimension of $X.$
Therefore, $A$ belongs to the subcategory
$\dhf{\D}.$

Suppose now that $A\in \dhf{\D}.$
The object $A$ is a bounded complex of coherent sheaves.
Let us take locally free  bounded above resolution
$P^{\cdot}\stackrel{\sim}{\to} A$ and consider
a good truncation $\tau^{\ge -k}P^{\cdot}$ for sufficient large $k\gg 0$
which is clearly isomorphic to $A$ in $\D.$

Since $A\in \dhf{\D}$
for any closed point $t: x\hookrightarrow X$ the groups
$\Hom(A, t_*\O_x[i])$ are zero for $|i|\gg 0$.
This means that
for sufficiently large $k\gg 0$ the truncation $\tau^{\ge -k}P^{\cdot}$
is a complex of locally free sheaves at the point $x,$ and, hence, in some neighborhood of $x.$
The scheme $X$ is quasicompact. This implies that there is a common
sufficiently large $k$ such that the truncation $\tau^{\ge -k}P^{\cdot}$
is a complex of locally free sheaves everywhere on $X,$ i.e. $A$ is perfect.
\end{dok}

\begin{corollary}\label{semde}
Let $X$ satisfies (ELF). Suppose that the
bounded derived category of coherent sheaves on $X$ has a
semiorthogonal decomposition
$
\db{\coh(X)}=\left\langle{\D_1}, \dots, {\D_n}\right\rangle.
$
Then the triangulated category of singularities
$\dsing{X}$ has the semiothogonal decomposition of the form
$$
\dsing{X}=\left\langle\ds{(\D_1)}, \dots, \ds{(\D_n)}\right\rangle.
$$
\end{corollary}
\begin{example}\label{exa}{\rm
Let $X$ be a projective space bundle $\PP(\E),$ where $\E$ is a vector bundle
of rank $r$ on a scheme $Y.$ It can be shown that the bounded derived category
$\db{\coh(X)}$ has a semiorthogonal decomposition of the form
$$
\db{\coh(X)}=\left\langle p^*\db{\coh(Y)}\otimes\O(i), \dots, p^*\db{\coh(Y)}\otimes\O(i+r-1)\right\rangle,
$$
for any $i\in\ZZ,$ where $p: X\lto Y$ is the projection.
This statement for a smooth
base can be found in \cite{Or1} and the proof works for any base.
By Corollary \ref{semde} the subcategory of perfect complexes
$\perf{X}$ also has the semiorthogonal decomposition
of the same form
$$
\perf{X}\cong
\bigl\langle p^*\perf{Y}\otimes \O(i), \dots,
 p^* \perf{Y\mathcal{}}\otimes\O(i+r-1)\bigl\rangle.
$$
And, finally, we obtain the semiorthogonal decomposition for the triangulated category
of singularities of $X$ in the terms of the triangulated category of singularities of $Y.$
}
\end{example}

Triangulated categories of singularities of $X$ have  an additional good properties
in  case of scheme is Gorenstein.
Recall that a local noetherian ring $A$ is called Gorenstein if
$A$ as a module over itself has a finite injective resolution.
It can be shown that if $A$ is Gorenstein than $A$ is a
dualizing complex for itself (see \cite{Ha}).
This means that $A$ has  finite injective dimension and
the natural map
$$
M\lto \bR\Hom^{\cdot}(\bR\Hom^{\cdot}(M, A), A)
$$
is an isomorphism for any coherent $A$\!-module $M$ and,
as a consequence, for any object from $\db{\coh(\Spec(A))}.$

\begin{definition} A scheme $X$ is  Gorenstein if
all of its local rings are Gorenstein local rings.
\end{definition}
If $X$ is Gorenstein and has finite dimension, then
$\O_X$ is a dualizing complex for $X,$ i.e. it has finite injective
dimension as the quasi-coherent sheaf and the natural map
$$
\F\lto \rhom^{\cdot}(\rhom^{\cdot}(\F, \O_X), \O_X)
$$
is an isomorphism for any coherent sheaf $\F.$
In particular, there is
an integer $n_0$ such that
$\ext^{i}(\F, \O_X)=0$ for each quasi-coherent sheaf
$\F$ and all $i>n_0.$

The following statements and their proofs can be found in \cite{Tr}.
\begin{lemma-definition}\label{rres}
Let X satisfy (ELF) and be Gorenstein. We say that
a coherent sheaf $\F$ is Cohen-Macaulay if the following  equivalent
conditions hold.
\begin{itemize}
\item[1)] The sheaves $\ext^{i}(\F, \O_X)$ are trivial for all
$i>0$.
\item[2)] There is a right locally free resolution
$
0\lto\F\lto\{Q^0\lto Q^1 \lto Q^3\lto \cdots \}.
$
\end{itemize}
\end{lemma-definition}
\begin{lemma}\label{locfr}
Let X satisfy (ELF) and be Gorenstein.
Let $\F$ be a Cohen Macaulay sheaf which is perfect as a complex.
Then $\F$ is locally free.
\end{lemma}
\begin{proposition}\label{fint}
Let $X$ satisfy (ELF) and be Gorenstein.
Then any object $A\in\dsing{X}$ is isomorphic
to the image of a Cohen-Macaulay sheaf $\F.$
\end{proposition}

\section{Reduction of dimension}

Let $S$ be a separated regular noetherian scheme of finite
Krull dimension. Let $\E$ be a vector bundle on $S$ of rank $r$ and
let $s \in H^0(S, \E)$
be a  section. Denote  by $X\subset S$ the zero subscheme of $s.$
We will assume that the section $s$ is regular, i.e. the codimension of
the subscheme $X$ coincides with the rank $r\footnote{By definition,
regularity means that the Koszul complex constructed
with $s$ is exact, but since $S$ is regular it is equivalent to
the codimension of the zero subscheme being the right one.}.$
In particular, the restriction of the bundle $\E$ on $X$ coincides
with the normal bundle to $X$ in $S$ which will be denoted
by $\N_{X/S}$

Consider the associated projective space bundles
$\PP(\E^{\vee})$ and $\PP(\N_{X/S}^{\vee}),$
where $\E^{\vee}$ and $\N_{X/S}^{\vee}$ are the dual vector bundles.
Denote these schemes
by $S'$ and $Z$ respectively, and the projections onto  $S$ and $X$ by $q$ and
$p$ respectively.
There are  canonical line bundles on $S'$ and $Z$ ,
denoted by $\O_{\E}(1)$ and $\O_{\N}(1),$ respectively, and the canonical surjections
\begin{equation}\label{twseq}
q^*\E\lto\O_{\E}(1),\qquad p^*\N_{X/S}\lto\O_{\N}(1).
\end{equation}
The section $s$ induces a section
$w\in H^0( S', \O_{\E}(1)).$
Denote by $Y$ the zero divisor on $S'$ defined by the section $w.$
The natural closed embedding of $Z =\PP(\N_{X/Y})$ into $S'$ goes through
$Y.$ Consider the closed embedding $i: Z\hookrightarrow Y.$
The kernel of the latter map of (\ref{twseq}) coincides with the normal bundle
to $N_{Z/Y},$ i.e. there is an exact sequence of sheaves on $Z$
\begin{equation}\label{Eul}
0\lto\N_{Z/Y}\lto p^*\N_{X/S}\lto\O_{\N}(1)\lto 0.
\end{equation}
%For an object $A$ from the derived category of coherent sheaves on $S'$ or on $Z$
%we often write $A(n)$ instead $A\otimes \O(n).$
All schemes defined above can be included in the following commutative diagram.
\begin{equation}\label{diag1}
\xymatrix{ Z=\PP(\N_{X/S}^{\vee}) \ar[dd]_{p}
\ar@{^{(}->}[r]^(.7){i}  & Y
 \ar[rdd]^{\pi}\ar@{^{(}->}[r]^(.35){u} & S'=\PP(\E^{\vee}) \ar[dd]^{q}
\\
\\
X \ar@{^{(}->}[rr]^{j} &  & S}
\end{equation}
Consider the composition functor $\bR i_* p^*:
\db{\coh(X)}\to \db{\coh(Y)}$ and denote it by $\Phi_{Z}.$
The aim of this section is to prove the following theorem.
\begin{theorem}\label{main1}
Let $S, S', X, Y,$ and $Z$ be as above. Then the functor
$\Phi_{Z}:\db{\coh(X)}\to \db{\coh(Y)}$ defined by the formula
$
\Phi_{Z}(\cdot)=\bR i_* p^*(\cdot)
$
induces a functor
$$
\ove{\Phi}_{Z}: \dsing{X}\lto \dsing{Y}
$$
which is an equivalence of triangulated categories.
\end{theorem}
We give two proofs of this theorem. Both of them use the following proposition.
\begin{proposition}\label{fff}
The functor $\Phi_{Z}=\bR i_* p^*:\db{\coh(X)}\to \db{\coh(Y)}$
 is fully faithful.
\end{proposition}
\begin{dok}
First, note that the functor $\Phi_{Z}=\bR i_{*} p^*$ has a right
adjoint functor which we denote by $\Phi_{Z *}.$
It can be represented as a composition $\bR p_{*} i^{\flat},$
where $i^{\flat}$ is right adjoint to $\bR i_{*}$
and has the form $\bL i^*(\cdot \otimes \omega_{Z/Y})[-r+1],$
where $\omega_{Z/Y}\cong \Lambda^{r-1}\N_{Z/Y}$
(see, for example, \cite{Ha} III, Cor.7.3).

Second, one can see that  the functor $p^*:\db{\coh(X)}\lto
\db{\coh(Z)}$ is fully faithful, because
$\bR p_{*} \O_{Z}\cong\O_{X}$ and by the projection
formula we obtain  isomorphisms
$$
\bR p_{*} p^*(A)\cong A \otimes \bR p_{*} \O_{Z}\cong
A\otimes \O_{X}\cong A
$$
for every $A\in \db{\coh(X)}.$

Now consider the canonical transformation of the functors
$\id \to \Phi_{{Z} *}\Phi_{Z}.$
To prove the proposition we have to show that this transformation
is an isomorphism.
This means we should check that for any object $A\in \db{\coh(X)}$
an object $C_A$ coming from an exact triangle
\begin{equation}\label{agtr}
A\lto \Phi_{{Z} *}\Phi_{Z} A\lto C_A
\end{equation}
is isomorphic to the zero object. All objects $A$ for which $C_A=0$ form
a triangulated subcategory of $\db{\coh(X)}.$
Since the minimal triangulated subcategory which contains all sheaves coincides
with whole $\db{\coh(X)}$ it is enough to show that $C_A=0$ when $A=\G$ is a sheaf.

Let $A=\G$ be a sheaf. Since the functor $p^*$ is fully faithful, the triangle (\ref{agtr})
is the image of the triangle
\begin{equation}\label{agtr1}
p^* \G\lto i^{\flat} \bR i_{*} p^* \G \lto B_{\G}
\end{equation}
under the functor $\bR p^*$, where the first morphism
is the canonical map induced by the natural transformation $\id\to i^{\flat}
\bR i_{*}$.
We have  the following sequence of isomorphisms
\begin{multline}
\bR i_* i^{\flat}\bR i_* p^*\G\cong
\bR i_* (\bL i^* \bR i_* p^*\G\otimes \omega_{Z/Y})[-r+1])\cong
\bR i_* (p^*\G\otimes\omega_{Z/Y})\otimes_{\O_Y} \O_Z [-r+1]\cong\\
\bR i_* (p^*\G\otimes\omega_{Z/Y})\stackrel{\bL}{\otimes}_{\O_Z}
(\O_Z\stackrel{\bL}{\otimes}_{\O_Y}\O_Z) [-r+1]).
\end{multline}

It is known (see \cite{Il}VII, 2.5) that
$
\tor_{k}^{\O_Y}(\O_Z, \O_Z)\cong
\Lambda^{k}\N_{Z/Y}^{\vee}.
$

Hence, the object $\bR i_* i^{\flat}\bR i_* p^*\G$ is a complex
the (r-k-1)-th cohomology of which is isomorphic to
$i_* (p^*\G\otimes\omega_{Z/Y}\otimes\Lambda^{k}\N_{Z/Y}^{\vee}).$
Therefore, the object $i^{\flat}\bR i_* p^*\G$ is a complex on $Z$
the $k$-th cohomology of which is isomorphic to
$p^*\G\otimes\Lambda^k\N_{Z/Y}$ for $k=0, \dots, r-1.$
This implies that the object $B_{\G}$ has  nontrivial
cohomologies only for $k=1, \dots, r-1,$ and  $H^k(B_{\G})$ is isomorphic to
$p^*\G\otimes\Lambda^k\N_{Z/Y}.$

Now, applying the functor $\bR p_*$ to the exact sequence
(\ref{Eul}) we obtain that $\bR p_{*} \N_{Z/Y}\cong 0.$ Moreover,
it can be easily checked that
\begin{equation}\label{Bott}
\bR p_* \Lambda^k \N_{Z/Y}\cong 0 \qquad\text{ for}\qquad k>0.
\end{equation}
(This is a relative analogue of the fact that on a projective space
$H^i(\Omega^k(-k))=0$ for all $i$ and all $k>0,$
which is a particular case of the Bott formula.)
The equalities (\ref{Bott}) implies that
$$
\bR p_{*} H^k(B_{\G})\cong 0
\quad
\text{for all}
\quad
k.
$$
Hence, $\bR p_* B_{\G}=0.$ Thus, the functor $\Phi_{Z}=\bR
i_{*} p^*$ is fully faithful.
\end{dok}

It can be useful to note that actually we have proven the following
more general proposition.
\begin{proposition} Let $p: Z=\PP(\N^{\vee}) \to X$ be the
projective space bundle associated to a vector bundle $\N$
on a scheme $X.$ Suppose that there is a regular closed embedding
$i: Z\hookrightarrow Y$ such that the normal
bundle $\N_{Z/Y}$ coincides with the kernel of the canonical map
$p^*\N\to \O_{\N}(1).$ Then the composition functor
$
\bR i_* p^* :\db{\coh(X)}\lto \db{\coh(Y)}
$
is fully faithful.
\end{proposition}

\begin{corollary}\label{ff}
The functor $\Phi_{Z}:\db{\coh(X)}\to \db{\coh(Y)}$
induces a functor
$$\ove{\Phi}_{Z}: \dsing{X}\lto \dsing{Y}$$
which is fully faithful too.
\end{corollary}
\begin{dok}
First, the functors $p^*$ and $i^{\flat}=\bL i^*(\cdot\otimes
\O(Z))[-1]$
take  perfect complexes to  perfect complexes as functors
of inverse images.
Second, the functors of direct images $\bR i_{*}$ and $\bR p_{*}$
also preserve perfect complexes, because both
 morphisms $i$ and $p$ have finite Tor-dimension.
Thus, by Lemma \ref{adjqu} we get a functor
$\ove{\Phi}_{Z}:\dsing{X}\to \dsing{Y}$ which is fully faithful.
\end{dok}

We  constructed the functor
$\ove{\Phi}_{Z}: \dsing{X}\lto \dsing{Y}$
and showed that this functor is fully faithful.

Now we give the first proof of Theorem \ref{main1}.
We will show  that  any object $F \in \dsing{Y}$ satisfying
condition $\ove{\Phi}_{Z*} F=0$ is isomorphic to the zero object.
The property for $\ove{\Phi}_{Z}$  to be an equivalence
is formally deduced from this fact.

The proof is  based on the following two simple lemmas.
\begin{lemma}\label{locdiv}
Let $i: Z\hookrightarrow Y$ be a closed embedding.
Let $\F$ be a coherent sheaf on $Y$ such that its restriction to
the complement $U=Y\setminus Z$ is locally free and $\bL i^* \F$
is isomorphic to a locally free sheaf on $Z$. Then $\F$ is locally
free on $Y$.
\end{lemma}
\begin{dok}
To prove that $\F$ is locally free it is sufficient to show that
for any closed point $t:y\hookrightarrow Y$ we have the equalities
$
\Ext^{i}(\F, t_*\O_y)=0
$
for all $i>0$. The sheaf $\F$ is locally free on $U$. Hence, we only need to
consider the case  $y\in Z$. This means that $t=i \cdot t'$ where $t':
y\hookrightarrow Z$ is closed embedding. In this case
$$
\Ext^{i}_{Y}(\F, t_*\O_y)=\Hom^{i}_{Z}(\bL i^* \F, t'_*\O_y)=0
$$
for $i>0$, because $\bL i^* \F$ is isomorphic to a locally free sheaf on $Z$.
\end{dok}
\begin{lemma}\label{semor}
An object $B\in \db{\coh(Z)}$ is perfect if and only if  $\bR
p_{*}(B(n))$ are perfect objects in $\db{\coh(X)}$ for all $n\in
\ZZ.$
\end{lemma}
\begin{dok}
If  the objects $B$ is perfect then the objects
$\bR p_{*} (B(n))$ are perfect too.
The inverse statement  follows from the semiorthogonal decomposition for
$\db{\coh(Z)}$ which was described in Remark \ref{exa}
\end{dok}
\begin{remark}{\rm
The statement of  Lemma \ref{semor}
remains true
for  any smooth morphism $Z\to X.$}
\end{remark}
Using these two lemmas we can prove the following proposition.
\begin{proposition}\label{zerob}
Assume that an object $F\in \dsing{Y}$ satisfies the condition
$\ove{\Phi}_{Z*} F\cong 0$. Then ${F}\cong 0$ in $\dsing{Y}$.
\end{proposition}
\begin{dok}
At the beginning note that all schemes $X, Z,$ and $Y$
are Gorenstein as they are locally complete intersections in  regular schemes.
By Proposition \ref{fint} we can assume
that the object $F$ is represented by a Cohen-Macaulay sheaf
$\F,$ in particular $\ext^{i}(\F, \O_Y)=0$ for all $i\ne 0$.
Note that any such $\F$ is locally free on the complement
$Y\setminus \Sing(Y)$.
In addition, for such $\F$ we have $\bL i^* \F\cong i^*\F$ is a sheaf,
because $\F$ has a right locally free resolution.

Denote by $\L$ the relatively ample line bundle on $Y$ obtained by
the restriction of the line bundle $\O_{\E}(1).$ Consider the object
$G=\bL u^* \bR u_* \F,$
where $u:Y\hookrightarrow S'$ is a closed embedding of the divisor $Y.$
On the one hand, the object $G$ is a perfect complex, because it is an inverse image
of a bounded complex of coherent sheaves from the smooth variety $S'.$
On the other hand, the complex $G$ has two cohomologies:  $\F$ in the zero place
and $\F\otimes \L^{-1}$ in the $(-1)$\!-st place.
Thus the image of $\F\otimes \L^{-1}$ in $\dsing{Y}$ is isomorphic to the image of
$\F[-2].$
By the assumption, the object $\ove{\Phi}_{Z*} \F$
is zero, hence $\ove{\Phi}_{Z*} (\F\otimes \L^n)=0$ for all $n\in \ZZ.$

Denote by $\F'$ the sheaf $i^* \F.$
We checked that $\bR p_{*}(\F'(n))$ are  perfect as  complexes
on $X$ for all $n\in \ZZ.$
By Lemma \ref{semor}  the sheaf $\F'=i^*\F$ is  perfect as the
complex on $Z$.

On the other hand, the sheaf $\F'=i^*\F$ has a right
locally free resolution, i.e. it is also a Cohen-Macaulay sheaf on $Z$.
Lemma \ref{locfr} implies now that $\F'=i^* \F=\bL i^* \F$
is locally free on $Z$.
Therefore, by Lemma \ref{locdiv} the sheaf $\F$ is
locally free on whole $Y$.
This means that $F$ is isomorphic to the zero
object in $\dsing{Y}.$
\end{dok}

\noindent{\bf The first proof of Theorem \ref{main1}}\quad
By Corollary \ref{ff} we already know that the functor $\ove{\Phi}_{Z}$ is fully faithful.
Now let us check that the functor
$\ove{\Phi}_{Z*}$ is also fully faithful.
Take an object $B\in \dsing{Y}$ and consider the natural map
$ \ove{\Phi}_{Z}\ove{\Phi}_{Z*} B\to B$.
Denote
by $C_B$ its cone. Applying the functor $\ove{\Phi}_{Z*}$
to the obtained exact triangle and using fully faithfulness
of the functor $\ove{\Phi}_Z,$
we get that $\ove{\Phi}_{Z*}C_B\cong 0$. By Proposition
\ref{zerob} the object $C_B$ is isomorphic to the zero object. Hence,
the map  $ \ove{\Phi}_{Z}\ove{\Phi}_{Z*} B\to B$ is an isomorphism
and the functor
$\ove{\Phi}_{Z*}$ is  fully faithful.
\hfill$\Box$

\begin{remark}\rm
Theorem \ref{main1} can be generalized to the case of quotient stacks.
Suppose that a group scheme $G$ acts on $S$ such that $\E$ can be considered
as an equivariant vector bundle with an invariant section $s.$
In this case we can extend the $G$\!-action to all schemes
$X,Y,Z$ and $S'.$ Moreover, we get a functor between
equivariant categories $\bR i_* p^*:
\db{\coh^G(X)}\to \db{\coh^G(Y)}.$
This functor induces  a functor
$\dsing{\left[X/G\right]}\lto \dsing{\left[Y/G\right]}$
which is an equivalence of triangulated categories as well.
\end{remark}

The second proof of Theorem \ref{main1} is based on a
semiorthogonal decomposition of the derived category of coherent sheaves on $Y.$
This kind of semiorthogonal decomposition is a particular case of more general
situation which will be described in the paper \cite{Ku}.
\begin{proposition}\label{dual}
The category $\db{\coh{Y}}$ has  a semiorthogonal decomposition
of the form
$$
\db{\coh{Y}}=\langle \bR i_*p^*\db{\coh(X)},\; \bL\pi^*\db{\coh(S)}\otimes\L,\;
 \dots,\; \bL\pi^*\db{\coh(S)}\otimes\L^{r-1}
\rangle,
$$
where $\L$ is the restriction of $\O_{\E}(1)$ on $Y.$
\end{proposition}
\begin{dok}
First, by Proposition \ref{fff} the functor $\bR i_*p^*$ is fully faithful.
Second, considering the sequence
$$
0\lto \O_{\E}(-1)\lto \O_{S'}\lto u_* \O_Y\lto 0.
$$

Applying the functors $\bR q_* (\O_{\E}(-k)\otimes -)$ to it,
we obtain that
$$
\bR\pi_* \O_Y\cong \bR q_* \bR u_* O_Y\cong \O_S
\qquad
\text{and}
\qquad
\bR\pi_* \L^{-k}=0
\quad
\text{for all}
\quad
0<k\le {r-2}.
$$
This implies that the functor
$\bL\pi^*:\db{\coh(S)}\lto \db{\coh(Y)}$ is fully faithful and
the sequence of subcategories
$
\left(\bL\pi^*\db{\coh(S)}\otimes\L, \dots, \bL\pi^*\db{\coh(S)}\otimes\L^{r-1}\right)
$
is semiorthogonal.

Third, let $A\in \db{\coh(X)}$ and $B\in\db{\coh(S)}.$
We have a sequence of isomorphisms
\begin{multline*}
\Hom(\bL\pi^* B\otimes \L^k,\; \bR i_*p^* A)\cong
\Hom(\bL i^*(\bL\pi^* B\otimes \L^k),\; p^* A)\cong
\Hom(p^*\bL j^* B,\; p^* A(-k))\cong\\
\Hom(\bL j^* B,\; A\otimes\bR p_* \O_{\N}(-k))=0.
\end{multline*}
And the last equality holds because $\bR p_*\O_{\N}(-k)=0$ for $0<k\le r-1.$
Thus we obtain the semiorthogonal sequence of subcategories
\begin{equation}\label{semi}
\left(
\bR i_*p^*\db{\coh(X)},\; \bL\pi^*\db{\coh(S)}\otimes\L,\;\dots,\;\bL\pi^*\db{\coh(S)}\otimes\L^{r-1}
\right).
\end{equation}

Finally, we should show that this sequence is full. Denote by $\C$ the subcategory of $\db{\coh(Y)}$
which is generated by the sequence (\ref{semi}). As $\C$ is admissible there is a semiorthogonal
decomposition $\D=\langle\C^{\perp}, \C\rangle.$
To show that $\C^{\perp}$ is trivial it is sufficient to check that all
 structure sheaves of closed points belong to $\C.$
Consider a closed point $y\in Y.$ Denote by $s$
the image $\pi(y)\in S.$
The fiber
over $s$ is isomorphic to the projective space $\PP^{r-2}$
if $s\in S\setminus X$ and it is isomorphic to $\PP^{r-1}$
if $s\in X.$
In the first case, since $\pi$ is flat over $s$ we have
$\O_{Y_s}\otimes \L^k=\bL\pi^*\O_s\otimes\L^k.$
These sheaves for $k=1,...,r-1$ belongs to $\C$
and they form an exceptional collection on
the projective space $Y_s=\PP^{r-2}.$
Therefore, the skyscraper sheaf $\O_y$ belongs to $\C$ too.

In the second case, we know that $\bR i_* p^* \O_s=\O_{Y_s}$ belongs to $\C$
and the objects
$\bL \pi^* \O_s\otimes \L^k=\bL j^* \O_{Y_s}\otimes\L^k$ also belong
to $\C$ for $k=1,\dots, r-1.$
The object $\bL \pi^* \O_s\otimes \L$ is a complex with two cohomologies
$\O_{Y_s}$ and $\O_{Y_s}\otimes\L.$ Since $\O_{Y_s}$ belongs to $\C$ then
$\O_{Y_s}\otimes\L\in \C$ too. Iterating this procedure we obtain that
$\O_{Y_s}\otimes\L^k\in\C$ for all $k=0,\dots, r-1.$ Since $Y_s=\PP^{r-1}$
 by the same reason as above we get that $\O_y\in\C$ for all $y\in Y.$
Thus our semiorthogonal sequence is full.
\end{dok}

\noindent{\bf The second proof of Theorem \ref{main1}}\quad
By Corollary \ref{semde} a semiorthogonal decomposition of the category
$\db{\coh(X)}$ induces the semiorthognal decomposition of $\dsing{Y}.$
Since $S$ is smooth $\dsing{S}$ is trivial. Therefore, we get that the category
$\dsing{Y}$ is equivalent to $\dsing{X}.$
\hfill$\Box$

\section{Application to D-branes in Landau-Ginzburg models.}

By a Landau-Ginzburg model we mean the following data: a smooth
variety $X$ with a symplectic K\"ahler form $\omega$
 and a regular nonconstant function $W$ on $X$
which is considered as a flat map $W: X\lto \AA^1$ and
which should be a symplectic fibration.
The function $W$ is called superpotential.
Note that for the
definition of D-branes of type B a symplectic form
is not needed.

A mathematical definition of the categories of D-branes of type B for  affine
Landau-Ginzburg models was proposed by M.Kontsevich
(see \cite{KL,Tr}). Roughly speaking, he
suggests that  superpotential $W$ deforms complexes of coherent
sheaves to "twisted" complexes, i.e the composition of
differentials  is no longer zero, but  is equal to multiplication
by $W.$ This "twisting" also breaks $\ZZ$\!-grading down to
$\ZZ/2$\!-grading. The equivalence of this definition with the
physics notion of B-branes in LG models was verified in  the paper
\cite{KL} in  the case of the usual quadratic superpotential
 and physical arguments were given
supporting Kontsevich's proposal for a general superpotential.

It was proved in the paper \cite{Tr} (Cor. 3.10)  that  the category of B-branes on
a smooth affine $X$ with a superpotential $W$ is equivalent to the product
$\mathop{\prod}\limits_{{\lambda}\in \AA^1} \dsing{X_{\lambda}},$
where $X_{\lambda}$ is the fiber over $\lambda\in\AA^1,$ and this
product is finite. For non-affine $X$ the category
$\mathop{\prod}\limits_{{\lambda}\in \AA^1} \dsing{X_{\lambda}}$ can be considered as a definition
of the category of D-branes of type B.

Let $S$ be a smooth quasiprojective variety
and let $f, g\in H^0(S, \O_S)$
be  two regular functions.
% such that
%the three functions $f,g$ and the constant $\b1$ would be lineary independent.
Suppose that the zero divisor $D\subset S$ defined by the function $g$
is smooth and the restriction of $f$ on $D$ is not constant.
We can consider $D$ as
a Landau-Ginzburg model with
superpotential
$f_{D}: D\lto \AA^1.$
Denote by $D_{\lambda}$ the fiber of the map $f_{D}: D\to\AA^1$
over a point $\lambda\in \AA^1.$

Another Landau-Ginzburg model
is given
by the smooth variety $T=S\times \AA^1$
and the superpotential $W: T \to \AA^1$
defined by the formula $W=f+xg,$ where $x$ is a coordinate on $\AA^1.$
Denote by $T_{\lambda}$ the fiber of $W$ over the point $\lambda.$
The natural closed embedding of $Z_{\lambda}=D_{\lambda}\times \AA^1$
into $T$ induces
the closed embedding $i_{\lambda}: D_{\lambda}\to T_{\lambda}.$
Consider the functor $$\Phi_{Z_{\lambda}}=\bR i_{\lambda*} p_{\lambda}^*:
\db{\coh(D_{\lambda})}\to \db{\coh(T_{\lambda})},$$ where $p_{\lambda}: Z_{\lambda}\to D_{\lambda}$
is the natural projection.
%Denote these functor by $\Phi_{Z_{\lambda}}.$
%
%
%
%
\begin{theorem}\label{sequn}
The functor $\Phi_{Z_{\lambda}}:\db{\coh(D_{\lambda})}\to \db{\coh(T_{\lambda})}$
induces a functor $\ove{\Phi}_{Z_{\lambda}}: \dsing{D_{\lambda}}\stackrel{\sim}{\lto} \dsing{T_{\lambda}}$
which is an equivalence of triangulated categories.
\end{theorem}
\begin{dok}
This is particular case of Theorem \ref{main1}.
To apply it we consider a trivial two dimensional vector bundle
$\E=\O_S^{\oplus 2}$ on $S$ and a section $s_{\lambda}\in H^0(S, \E)$
which is given by the two functions $g, f-\lambda\b1.$ Thus $D_{\lambda}$
is the zero subvariety of $s_{\lambda}$ and it is the analogue of the scheme $X$ from Theorem \ref{main1}.
Applying the construction preceding that theorem
we obtain varieties $Z=X\times \PP^1$ and $Y\subset S\times\PP^1$
which are compactifications of $Z_{\lambda}$ and $T_{\lambda}$ respectively.
By Theorem \ref{main1} triangulated categories of singularities of $X=D_{\lambda}$ and $Y$
are equivalent.
The complement of $T_{\lambda}$ in $Y$ is a relative ample divisor
which does not meet $\Sing(Y)=\Sing(T_{\lambda}).$
Hence, by Proposition \ref{locpr} $Y$ and $T_{\lambda}$ also have equivalent triangulated categories of singularities.
\end{dok}
\begin{corollary}\label{LG} Let $S$ be a smooth quasi-projective variety. Then
the Landau-Ginzburg models $f_D: D\to\AA^1$ and $W: S\times \AA^1\to \AA^1,$
where $W=f+xg$ and $D:=\{g=0\}\subset S$ is smooth,
have equivalent categories of D-branes of type B.
\end{corollary}
%\vspace{-0cm}

\end{document}